\documentclass{article}
\usepackage[accepted]{icml2009}
\usepackage{epsfig,amsmath,amssymb,epsfig}

\icmltitlerunning{Orbit-Products \& Gaussian BP}

\newcommand{\tr}{\mathrm{tr}}
\newcommand{\qed}{$\diamond$}
\newcommand{\bp}{\mathrm{bp}}
\newcommand{\diff}{\setminus}

\newcommand{\ignore}[1]{}

\newtheorem{theorem}{Theorem}
\newtheorem{lemma}{Lemma}
\newtheorem{corollary}{Corollary}

\begin{document}

\twocolumn[
\icmltitle{Orbit-Product Representation and Correction of\\ Gaussian Belief
Propagation}
\icmlauthor{Jason K. Johnson$^a$}{jasonj@lanl.gov}
\icmlauthor{Vladimir Y. Chernyak$^{b,a}$}{chernyak@chem.wayne.edu}
\icmlauthor{Michael Chertkov$^a$}{chertkov@lanl.gov}
\icmladdress{$^a$Theoretical Division and Center for Nonlinear Studies, LANL,
             Los Alamos, NM 87545, USA}
\vspace{-.3cm}
\icmladdress{$^b$Department of Chemistry, Wayne State University,
             Detroit, MI 48202, USA}
\vskip 0.3in
]

\begin{abstract} We present a new view of Gaussian belief propagation (GaBP)
based on a representation of the determinant as a product over orbits of a
graph.  We show that the GaBP determinant estimate captures totally backtracking
orbits of the graph and consider how to correct this estimate.  We show that the
missing orbits may be grouped into equivalence classes corresponding to
backtrackless orbits and the contribution of each equivalence class is easily
determined from the GaBP solution. Furthermore, we demonstrate that this
multiplicative correction factor can be interpreted as the determinant of a
backtrackless adjacency matrix of the graph with edge weights based on GaBP.
Finally, an efficient method is proposed to compute a truncated correction
factor including all backtrackless orbits up to a specified length.
\end{abstract}

\section{Introduction}

Belief Propagation is a widely used method for inference in graphical models. We
study this algorithm in the context of Gaussian graphical models.  There have
been several studies of Gaussian belief propagation (GaBP)
\cite{Weiss01,Rusmevichientong01,Plarre04} as well as numerous applications
\cite{Moallemi06,Bickson08a,Bickson08b}. The best known sufficient condition for
its convergence is the walk-summable condition \cite{Johnson06,Malioutov06} (see
also \cite{Cseke08,Moallemi09}), which also provides new insights into the
algorithm by interpreting it as computing weighted sums of walks (walk-sums)
within the graph. Our aim in this present paper is to extend this
graphical/combinatorial view of GaBP to include estimation of the determinant
(partition function) of the Gaussian graphical model.  This work is also
inspired by the loop-series correction method for belief propagation
\cite{Chertkov06} that was recently extended to Gaussian graphical models
\cite{Chernyak08}.

Our present study leads to a new perspective on GaBP having close ties to
graphical zeta functions \cite{Stark96}. We find that for walk-summable models
the determinant may be represented as a product over all orbits (cyclic walks)
of the graph.  The estimate of the determinant provided by GaBP only captures
totally backtracking orbits, which can be embedded as orbits in the
computation tree (universal cover) of the graph.  The missing orbits may then be
grouped into equivalence classes corresponding to backtrackless orbits.  The
orbit-product over each such equivalence class may be simply computed from the
solution of GaBP. Also, the product over all backtrackless orbits may be
interpreted as the determinant of a backtrackless adjacency matrix of the graph
with appropriately defined edge weights based on the GaBP solution. Finally, we
propose a simple, efficient method to compute truncated orbit-products including
all orbits up to some specified length and provide an error-bound on the
resulting estimates. In certain classes of graphs (e.g., grids), this leads to
an efficient method with complexity linear in the number of nodes and the
required precision of the determinant estimate.

This paper differs fundamentally from \cite{Chernyak08} in that we rely heavily
on the walk-summable property to develop multiplicative expansions using
infinitely many orbits of the graph, whereas \cite{Chernyak08} develops additive
expansions over a finite number of ``generalized loops'' (which may be
disconnected) using methods of Grassman calculus.  Our present approach leads
naturally to approximation methods (with accuracy guarantees in walk-summable
models) based on truncated orbit-products.

\section{Preliminaries}

\subsection{Walks and Orbits of a Graph}

Let $G$ be a graph on vertices (nodes) $V = \{1,\dots,n\}$ with undirected edges
$\{i,j\} \in G$.  We may also treat each undirected edge $\{i,j\}$ as a
symmetric pair of directed edges $(ij)$ and $(ji)$. A \emph{walk} $w$ is a
sequence of adjacent vertices $(w_0 \dots w_L)$ ($w_t  \in V$ for $t=0,\dots,L$
and $\{w_t,w_{t+1}\} \in G$ for $t = 1,\dots,L-1$) where $|w| \triangleq L$ is
the \emph{length} of the walk. A walk may be equivalently specified as a
sequence of \emph{steps} $w = ((w_0 w_1)(w_1 w_2)\dots(w_{L-1} w_L))$ such that
each
step is a (directed) edge of the graph that ends where the following step
begins.  A walk may visit the same node or cross the same edge multiple times
and may also \emph{backtrack}---that is, it may step back to the preceding
vertex.  A walk is \emph{closed} if it begins and ends at the same node $w_0 =
w_L$.  A closed walk is \emph{primitive} if it is not a multiple of some shorter
walk (e.g., the walk $(1231)$ is primitive but $(1231231)$ is not).  We define
an \emph{orbit} $\ell = [w]$ to be an equivalence class of closed primitive
walks, where two walks are considered equivalent $[w] = [w']$ if one is a
cyclic shift of the other (i.e., if $w_t = w'_{t+s(\mathrm{mod} \, \ell)}$ for
some $s$ and $t=1,\dots,L$).  Hence, there is a one-to-one correspondence
between orbits and (non-terminating) cyclic walks.

The following classification of walks and orbits plays an essential role in our
analysis: A walk (or orbit) is said to be \emph{reducible} if it contains a
\emph{backtracking pair} of consecutive steps $\dots(ij)(ji)\dots$, otherwise
the walk is \emph{irreducible}.  By repeatedly deleting backtracking pairs until
none remain one obtains the (unique) \emph{irreducible core} $\gamma =
\Gamma(w)$ of the walk $w$. For closed walks it may happen that $\gamma =
\emptyset$,  where $\emptyset \triangleq ()$ denotes the trivial (empty) walk.
We then say that the walk is \emph{totally reducible}. We say that a walk is
\emph{non-trivial} if it is not totally reducible. Totally reducible walks have
been called \emph{backtracking} \cite{Malioutov06},  although \emph{totally
backtracking} is perhaps a better description. Irreducible walks have been
called \emph{backtrackless} (or \emph{non-backtracking}) elsewhere in the
literature.

\emph{Notation.} $\mathcal{L}$ denotes the set of all orbits of $G$,
$\Gamma(\mathcal{L})$ denotes the irreducible (backtrackless) orbits, and we
partition $\mathcal{L}$ into (disjoint) equivalence classes $\mathcal{L}_\gamma
\triangleq \{ \ell \in \mathcal{L} | \Gamma(\ell) = \gamma \}$ for $\gamma \in
\Gamma(\mathcal{L})$.   In particular, $\mathcal{L}_\emptyset$ denotes the class
of totally backtracking orbits, which plays a special role in our interpretation
of GaBP.  We will use $\ell$ to denote a generic orbit and reserve $\gamma$ to
denote irreducible orbits. 

\emph{Example.} Orbits $[1231]$, $[1231451]$, $[123421561]$ are
backtrackless; $[1234321]$, $[1232421]$, $[1231321]$ are totally backtracking; $[123241]$ 
is both reducible and non-trivial (neither backtrackless nor totally backtracking).

\subsection{Gaussian Belief Propagation}

A \emph{Gaussian graphical model} is a probability distribution \begin{equation}
p(x) = \mathcal{Z}^{-1} \exp\left\{ -\tfrac{1}{2} x^T J x + h^T x \right\}
\end{equation} of random variables $x \in \mathbb{R}^n$ where $J$ is a sparse,
symmetric, positive-definite matrix.  The fill-pattern of $J$ defines a graph
$G$ with vertices $V = \{1,\dots,n\}$ and edges $(ij)$ for all $J_{ij} \neq 0$.
The partition function is defined by $\mathcal{Z}(h,J) \triangleq \int \exp\{
-\frac{1}{2} x^T J x + h^T x \} dx = \left[ (2\pi)^n \det J^{-1} \right]^{1/2}
e^{\tfrac{1}{2} h^T J^{-1} h}$ so as to normalize the distribution.  Given such
a model, we may compute the mean vector $\mu \triangleq \int p(x) x dx = J^{-1}
h$ and covariance matrix $K \triangleq \int p(x) (x-\mu)(x-\mu)^T dx = J^{-1}$.
This generally requires $\mathcal{O}(n^3)$ computation in dense graphs using
Gaussian elimination.  If $G$ is sparse and only certain elements of $K$ are
required (the diagonal and edge-wise covariances), then the complexity of
Gaussian elimination may be substantially reduced (e.g., $\mathcal{O}(n^{3/2})$
in planar graphs using nested dissection) but still generally has complexity
growing as $\mathcal{O}(w^3)$ in the \emph{tree-width} $w$ of the graph.

\emph{Gaussian belief propagation} (GaBP) is a simple, distributed, iterative
message-passing algorithm to estimate the marginal distribution $p(x_i)$ of each
variable, which is specified by its mean $\mu_i$ and variance $K_{ii}$. GaBP is
parameterized by a set of messages $m_{ij}(x_j) = e^{ \frac{1}{2} \alpha_{ij}
x_j^2 + \beta_{ij} x_j}$ defined on each directed edge $(ij)$ of the graph
($m_{ij}$ is regarded as a message being passed from $i$ to $j$). The GaBP
equations are:
\begin{displaymath}
m_{ij}(x_j) \propto \int \psi_i(x_i) \prod_{k \in \partial i \diff j}
m_{ki}(x_i)
\psi_{ij}(x_i,x_j) dx_i
\end{displaymath}
where $\psi_i = e^{-\tfrac{1}{2}J_{ii}x_i^2+h_i x_i}$, $\psi_{ij} =
e^{-J_{ij}x_ix_j}$ and $\partial i$ denotes the set of neighbors of $i$ in $G$. 
This reduces to the following rules for computing $(\alpha,\beta)$-messages:
\begin{eqnarray}
\alpha_{ij} &=&  J_{ij}^2 (J_{ii} - \alpha_{i \diff j})^{-1} \nonumber\\
\beta_{ij} &=&  -J_{ij} (J_{ii} - \alpha_{i \diff j})^{-1} (h_i + \beta_{i \diff
j}) \label{eq:gabp_eqs}
\end{eqnarray}
where $\alpha_{i \diff j} \triangleq \sum_{k \in \partial i \diff j}
\alpha_{ki}$ and $\beta_{i \diff j} \triangleq \sum_{k \in \partial i \diff j}
\beta_{ki}$.  These equations are solved by iteratively recomputing each message
from the other messages until convergence. The marginal distribution is then
estimated as $p^\bp(x_i) \propto \psi_i(x_i) \prod_{k \in \partial i}
m_{ki}(x_i)$, which gives
variance estimates $K^\bp_i = (J_{ii}-\sum_{k} \alpha_{ki})^{-1}$
and mean estimates $\mu^\bp_i = K^\bp_i (h_i+\sum_{k}
\beta_{ki})$. In trees, this method is equivalent to Gaussian elimination,
terminates after a finite number of steps and then provides the correct
marginals. In loopy graphs, it may be viewed as performing Gaussian elimination
in the computation tree (universal cover) of the graph
\cite{Plarre04,Malioutov06} (obtained by ``unrolling'' loops) and may therefore
fail to converge due to the
infinite extent of the computation tree. If it does converge, the mean estimates
are still correct but the variances are only approximate. We also obtain an
estimate of the pairwise covariance matrix on edges $\{i,j\} \in G$:
\vspace{-.15cm}
\begin{displaymath}
K^\bp_{(ij)} =
\left(
\begin{array}{cc}
J_{ii}-\alpha_{i \diff j} & J_{ij} \\
J_{ij} & J_{jj}-\alpha_{j \diff i}
\end{array}
\right)^{-1}
\end{displaymath}
In this paper, we are concerned with the GaBP estimate of the determinant $Z =
\det K = \det J^{-1}$ (which is closely linked to computation of the partition
function $\mathcal{Z}$).  We obtain an estimate of $Z$ from the GaBP solution
as:
\vspace{-.1cm}
\begin{equation} \label{eq:Z_bp_def}
Z^\bp = \prod_{i \in V} Z_i^\bp \prod_{\{i,j\} \in G} \frac{Z_{ij}^\bp}{Z_i^\bp
Z_j^\bp}
\end{equation}
where $Z_i^\bp = K^\bp_i$ and $Z_{ij}^\bp = \det K^\bp_{(ij)}$.   The motivation
for this form of estimate is that it becomes exact if $G$ is a tree. In loopy
graphs, there may generally be no stable solution to GaBP (or multiple unstable
solutions).  The main objective of this paper is to interpret the estimate
$Z^\bp$ in the context of walk-summable models (described below), for which
there is a well-defined stable solution, and to suggest methods to correct this
estimate. Note that the variance and covariance estimates (and hence the
determinant estimate $Z^\bp$) are independent of $h$ and the $\beta$-messages in
(\ref{eq:gabp_eqs}), they are determined solely by the $\alpha$-messages
determined by $J$.  Since GaBP correctly computes the means in walk-summable
models, we are mainly concerned with how to correct $Z^\bp$ (and hence its
derivatives, which correspond to the GaBP estimates of variances/covariances).

\paragraph{Walk-Sum Interpretation}

Our approach in this paper may be considered as as extension of the walk-sum
interpretation of GaBP \cite{Malioutov06}. Let $J$ be normalized to have
unit-diagonal, such that $J=I-R$ with $R$ having zeros along its diagonal. The
walk-sum idea is based on the series $K = (I-R)^{-1} = \sum_k R^k$, which
converges if $\rho(R)<1$ where $\rho(R)$ denotes the \emph{spectral radius} of
the matrix $R$ (the maximum modulus of the eigenvalues of $R$). This allows us
to interpret $K_{ij}$ as a sum over all walks in the graph $G$ which begin at
node $i$ and end at node $j$ where the weight of a walk is defined as $R^w =
\prod_{(ij) \in w} r_{ij}^{n_{ij}(w)}$ and $n_{ij}(w)$ is a count of how many
times step $(ij)$ occurs in the walk.  We write this \emph{walk-sum} as $K_{ij}
= \sum_{w: i \rightarrow j} R^w$.  However, in order for the walk-sum to be
well-defined, it must converge to the same value regardless of the order in
which we add the walks.  This is equivalent to requiring that it converges
absolutely.  Thus, we say that $R$ is \emph{walk-summable} if $\sum_{w: i
\rightarrow j} |R^w|$ converges for all $i,j \in V$.  This is equivalent to the
spectral condition that $\rho(|R|)<1$ where $|R| \triangleq (|r_{ij}|)$ is the
element-wise 
absolute-value matrix of $R$. A number of other equivalent or
sufficient conditions are given in \cite{Malioutov06}.

In walk-summable models it then holds that variances correspond to closed
walk-sums $K_{ii} = \sum_{w: i \rightarrow i} R^w$ and means correspond to a
(reweighted) walk-sum over all walks which end at a specific node $\mu_i =
\sum_{w: * \rightarrow i} h_* R^w$ (here $*$ denotes the arbitrary starting
point of the walk).  Moreover, we may interpret the GaBP message parameters
$(\alpha,\beta)$ as recursively computing walk-sums within the \emph{computation
tree} \cite{Malioutov06}.  This implies that GaBP converges in walk-summable
models and converges to the same ``walk-sum'' solution independent of the order
in which we update messages. This interpretation also shows that GaBP computes
the correct walk-sums for the means but only computes a subset of the closed
walks needed for the variances. Specifically, $K^\bp_i$ only includes totally
backtracking walks at node $i$. This is seen as a walk is totally backtracking
if and only if it can be embedded as a closed walk in the computation tree of
the graph and it is these closed walks of the computation tree which GaBP
captures in its variance estimates. 

\section{Orbit-Product Interpretation of Gaussian BP}

\subsection{Determinant $Z$ as Orbit-Product}

Let $Z(R) \triangleq \det(I-R)^{-1}$.  In walk-summable models, we may give this
determinant another graphical interpretation as a product over orbits of a
graph, one closely related to the so-called \emph{zeta function} of a graph
\cite{Stark96}.

\begin{theorem} If $\rho(|R|)<1$ then it holds that $Z(R) =
  \prod_\ell (1-R^\ell)^{-1} \triangleq \prod_\ell Z_\ell$ where the
  product is taken over all orbits of $G$ and $R^\ell = \prod_{(ij)
    \in \ell} r_{ij}^{n_{ij}(\ell)}$ where $n_{ij}(\ell)$ is the
  number of times step $(ij)$ occurs in orbit $\ell$.
\end{theorem}

\emph{Proof.} $\log \det (I-R)^{-1}
 = \tr \log (I-R)^{-1}
 = \tr \sum_k \frac{R^k}{k} \linebreak
 = \sum_{\mathrm{closed}\, w} \frac{R^w}{|w|}
 = \sum_{\mathrm{primitive}\, w} \sum_{m=1}^\infty \frac{(R^w)^m}{m |w|}
 = \sum_{\mathrm{primitive}\, w} \frac{1}{|w|} \log (1-R^w)^{-1}
 = \sum_{\mathrm{orbits}\, \ell} \log (1-R^\ell)^{-1}
 = \log \prod_\ell (1-R^\ell)^{-1}$.
We have used the identity $\log\det
A = \tr\log A$ and the series expansion
$\log (I-A)^{-1} = \sum_{k=1}^\infty \frac{A^k}{k}$.
Each closed walk is expressed as a multiple of a primitive
walk. Every primitive walk $w$ has exactly $|w|$ distinct cyclic shifts. \qed

We emphasize that $\rho(|R|)<1$ is necessary to insure that the the
orbit-products we consider are well-defined. \emph{This condition is assumed
throughout the remainder of the paper.}

\subsection{$Z^\bp$ as Totally Backtracking Orbits}

Totally backtracking walks play an important role in the walk-sum interpretation
of the GaBP variance estimates. We now derive an analogous interpretation of
$Z^\bp$ defined by (\ref{eq:Z_bp_def}):

\begin{theorem} $Z^\bp = \prod_{\ell \in \mathcal{L}_\emptyset} Z_\ell$
  where the product is taken over the set of totally backtracking orbits
  of $G$.
\end{theorem}

Although this result seems intuitive in view of prior work, its
proof is non-trivial involving arguments not used previously.  
To prove the theorem, we first summarize some useful lemmas.
Consider a block matrix $A = (A_{11} A_{12}; A_{21} A_{22})$.  
The \emph{Schur complement} of block $A_{11}$ is
$A_{22}^* \triangleq A_{22}-A_{21} A_{11}^{-1} A_{21}$. It holds
that $\det A = \det A_{11} \det A_{22}^*$ and $(A_{22}^*)^{-1} = (A^{-1})_{22}$.
Using these well-known identities, it follows:

\begin{lemma} Let $R = (R_{11} R_{12}; R_{21} R_{22})$ and $K = (I-R)^{-1} =
  (K_{11} K_{12}; K_{21} K_{22})$.  Then $\det K_{11} = \frac{Z(R)}{Z(R_{22})}$.
\end{lemma}
For walk-summable models, we then have
\begin{displaymath}
\det K_{11} = \frac{\prod_{\ell \in G} Z_\ell}{\prod_{\ell \in G_2} Z_\ell}
= \prod_{\ell \in G | \ell \,\mathrm{intersects}\, G_1} Z_\ell
\end{displaymath}
where the final orbit-product is taken over all orbits of $G$ which include any 
node of subgraph $G_1$ (corresponding to submatrix $R_{11}$). Next, using this
result
and the interpretation of GaBP as inference on the computation tree, we
are led to the following interpretation of the quantities $Z_i^\bp$ and 
$Z_{ij}^\bp$ appearing in (\ref{eq:Z_bp_def}). 
Let $T_i$ denote the computation tree of the
graph $G$ with one copy of node $i$ marked.  Let $T_{ij}$ denote the computation
tree with one copy of edge $\{i,j\} \in G$ marked. Then,

\begin{lemma} $Z_i^\bp = \prod_{\ell \in T_i | i \in \ell} Z_\ell$
  where the product is over all orbits of $T_i$ that include the
  marked node $i$.  $Z_{ij}^\bp = \prod_{\ell \in T_{ij} | i \in \ell
\,\mathrm{or}\, j \in \ell} Z_\ell$ where the product is
  over all orbits of $T_{ij}$ that include either endpoint of the
  marked edge $\{i,j\}$.
\end{lemma}

\ignore{\emph{Proof.} Let $R_T$ denote the weighted adjacency matrix of the
(infinite) computation tree and $K_T = (I-R_T)^{-1}$ its covariance
matrix.  Let node $i$ be marked in the computation tree. Then $Z_i^\bp
= \det (K_T)_{ii}$. Applying Lemma 2 and Theorem 1, we have:
\begin{displaymath}
Z_i^\bp
= \frac{Z(R_T)}{Z(R_{T \diff i})}
= \frac{\prod_{\ell \in T_i} Z_\ell}{\prod_{\ell \in T_i | i \not\in \ell}
Z_\ell}
= \prod_{\ell \in T_i | i \in \ell} Z_\ell
\end{displaymath}
A similar argument holds for $Z_{ij}^\bp$:
\begin{displaymath}
Z_{ij}^\bp
= \frac{\prod_{\ell \in T_{ij}} Z_\ell}
{\prod_{\ell \in T_{ij} | i \not\in \ell \wedge j\not\in \ell} Z_\ell}
= \prod_{\ell \in T_{ij} | i \in \ell \vee j \in \ell} Z_\ell \;\;\; \diamond
\end{displaymath}}

\emph{Proof of Theorem 2.} Using Lemma 2 and the correspondence between orbits
of the computation tree and totally backtracking orbits of $G$, we may expand
(\ref{eq:Z_bp_def}) to express  $Z^\bp$ entirely as a product over totally
backtracking orbits $Z^\bp = \prod_{\ell \in \mathcal{L}_\emptyset}
Z_\ell^{N_\ell}$ where $N_\ell$ is the count of how many times $\ell$
appears in the orbit-product---the number of times it appears in the numerator
of (\ref{eq:Z_bp_def}) minus the number of times in appears in the denominator.
It remains to show
that $N_\ell=1$ for each totally backtracking orbit.  This may be seen by
considering the subtree $T_\ell$ of the computation tree $T$ traced out by orbit
$\ell$. Let $v$ and $e$ respectively denote the number of nodes and edges of
$T_\ell$ (hence, $e = v-1$) and let $c$ denote the number of edges of $T$ with
exactly one endpoint in $T_\ell$.  First, we count how many powers of $Z_\ell$
appear in the orbit product $\prod_i Z_i^\bp$.  For each vertex $i \in T_\ell$
we may pick this as the marked node in the computation tree and this shows one
way that $\ell$ can be embedded in $T_i$ so as to include its marked node.
Thus, $v$ gives the total number of multiples of $Z_\ell$ in $\prod_i Z_i^\bp$.
Similarly, we could mark any edge $\{i,j\} \in T$ with one or both endpoints in
$T_\ell$ and this gives one way to embed $\ell$ into $T_{ij}$ so as to
intersect the marked edge.  Thus, the product $\prod_{ij} Z_{ij}^\bp$
contributes $e+c$ powers of $Z_\ell$.  Lastly, the product $\prod_{ij} Z_i^\bp
Z_j^\bp$ contains $2e+c$ powers of $Z_\ell$.  This represents the number of ways
we may pick a \emph{directed} edge $(ij)$ of $T$ such that at least one
endpoint is in $T_\ell$. The total count is then $N_\ell = v + (e+c) - (2e+c) =
v-e = 1$. \qed

Combining Theorems 1 and 2, we obtain the following orbit-product correction
to $Z^\bp$:
\begin{corollary} $Z = Z^\bp \times \prod_{\ell \not\in \mathcal{L}_\emptyset}
Z_\ell$.
\end{corollary}
This formula includes a correction for every missing orbit, that is,
for every non-trivial orbit.  This implies that $Z = Z^\bp$ for 
trees since \emph{all} orbits of trees are totally backtracking.

\subsection{$Z^\bp$ Error Bound}

One useful consequence of the orbit-product interpretation of $Z^\bp$ is that it
provides
a simple error bound on GaBP. Let $g$ denote the \emph{girth} of the graph $G$, 
defined as the length of the shortest cycle of $G$.  We note that the missing
orbits $\ell \not\in \mathcal{L}_\emptyset$ must all have length greater than or
equal to
$g$.  Then,
\begin{theorem}
$\frac{1}{n} \left|\log \frac{Z^\bp}{Z}\right| \le \frac{\rho(|R|)^g}{g
(1-\rho(|R|))}$.
\end{theorem}
\emph{Proof.}  We derive the chain of inequalities: 
$\left|\log \frac{Z}{Z^\bp} \right| 
\stackrel{\mathrm{(a)}}{=} \left| \sum_{\ell \not\in \mathcal{L}_\emptyset} \log
Z_\ell \right| 
\le \sum_{|\ell| \ge g} |\log Z_\ell| 
\stackrel{\mathrm{(b)}}{\le} \sum_{|\ell| \ge g} \log (1-|R|^\ell)^{-1} 
\stackrel{\mathrm{(c)}}{=} \tr \sum_{k=g}^\infty \frac{|R|^k}{k}
\le n \sum_{k=g}^\infty \frac{\rho^k}{k}
\le \frac{n \rho^g}{g} \sum_{k=0}^\infty \rho^k
=  \frac{n \rho^g}{g (1-\rho)}$. 
(a) Corollary 1.
(b) $|\log Z_\ell| 
= \left| \sum_{k=1}^\infty \frac{(R^\ell)^k}{k} \right| 
\le \sum_{k=1}^\infty \frac{(|R|^\ell)^k}{k}
= \log (1-|R|^\ell)^{-1}$. 
(c) The proof of Theorem 1 shows that $\tr \sum_{k \ge 1} \frac{R^k}{k} = 
\sum_{\ell} \log (1-R^\ell)^{-1}$.  Similarly,
$\tr \sum_{k \ge g} \frac{|R|^k}{k} 
= \sum_{|\ell| \ge g} \log (1-|R|^\ell)^{-1}$.
$\diamond$

This is consistent with the usual intuition that belief propagation is most
accurate in large girth graphs with weak interactions.

\section{Backtrackless Orbit Correction}

In this section we show that the set of orbits omitted in the GaBP estimate can
be grouped into equivalence classes corresponding to backtrackless
orbits and that the orbit-product over each such equivalence class is simply
computed with the aid of the GaBP solution:

\begin{theorem} $Z = Z^\bp \times \prod_{\gamma \neq \emptyset}  Z'_\gamma$
where the
product is over all \emph{backtrackless orbits} of $G$ and we define
\begin{displaymath}
Z'_\gamma = (1 - \prod_{(ij) \in \gamma} (r'_{ij})^{n_{ij}(\gamma)})^{-1}
\end{displaymath}
where $r'_{ij} \triangleq \frac{r_{ij}}{1 - \alpha_{i \diff j}}$ and
$\alpha_{i \diff j} = \sum_{k \in \partial i \diff j} \alpha_{ki}$ is computed
from the solution of GaBP.
\end{theorem}

In comparison to Corollary 1, here the correction factor is expressed as an
orbit-product over just the backtrackless orbits (whereas Corollary 1 uses a
separate correction for each non-trivial orbit).  However, all orbits are still
correctly accounted for because we modify the edge weights of the graph so as to
include a factor $(1-\alpha_{i \diff j})^{-1}$ (computed by GaBP) which serves
to ``factor in'' totally-backtracking excursions at each point along the
backtrackless orbit, thereby generating all non-trivial orbits.

The basic idea underlying this construction is depicted in Figure 1.  For each
backtrackless orbit $\gamma$ we define an associated \emph{computation graph}
$G_\gamma$ as follows. First, we start with a single directed cycle based on
$\gamma = [\gamma_1\gamma_2\cdots\gamma_L]$ (any duplicated nodes of the orbit
map to distinct nodes in this directed cycle).  Then, for each node $\gamma_k$
of this graph, we attach a copy of the computation tree $T_{\gamma_k \diff
\gamma_{k+1}}$, obtained by taking the full computation tree $T_{\gamma_k}$
rooted at node $\gamma_k$ and deleting the branch $(\gamma_k,\gamma_{k+1})$
incident to the root.  This construction is illustrated in Figure 1(a,b) for the
graph $G = K_4$ and orbit $\gamma = [(12)(23)(31)]$.  The cycle has ``one-way''
directed edges whereas each computation tree has ``two-way'' undirected edges.
This is understood to mean that walks are allowed to backtrack within the
computation tree but not within the cycle.  The importance of this graph is
based on the following lemma (the proof is omitted):

\begin{figure}
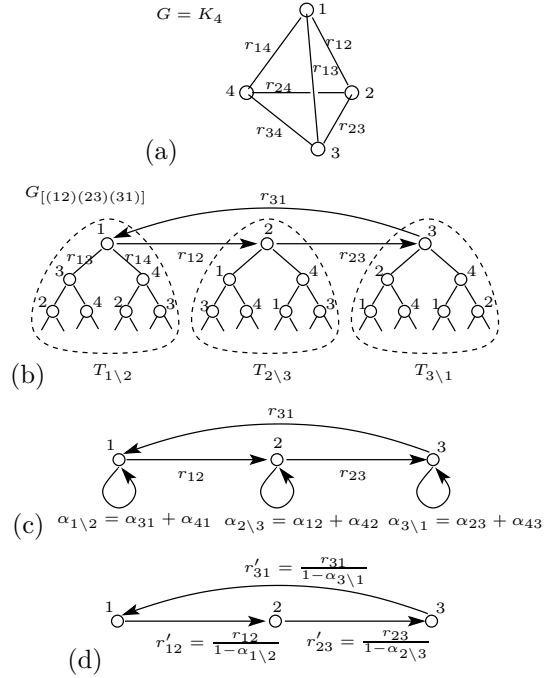

\centering
(a)\hspace{-.3cm}\input{K4.pstex_t}\\
\vspace{.1cm}
(b)\hspace{-.3cm}\input{loop_graph.pstex_t}\\
\vspace{.1cm}
(c) \input{loop_graph2.pstex_t}\\
\vspace{.1cm}
(d)\hspace{-.1cm}\input{loop_graph3.pstex_t}\\
\vspace{-.2cm}
\caption{Illustration of construction to combine equivalent orbits.  (a) The
graph $G=K_4$. (b) The computation graph $G_\gamma$ for $\gamma =
[(12)(23)(31)]$. (c) Finite graph with self-loops at each node to capture
totally backtracking walks. (d) Equivalent graph with modified edge weights to
capture totally backtracking walks.}
\vspace{-.3cm}
\end{figure}

\begin{lemma} Let $\gamma$ be a backtrackless orbit of $G$.  Then,
  there is a one-to-one correspondence between the class of orbits
  $\mathcal{L}_\gamma$ of $G$ and the non-trivial orbits of $G_\gamma$.
\end{lemma}

Next, we demonstrate how to compute all of the orbits within an equivalence
class as a simple determinant calculation based on the backtrackless orbit
$\gamma$ and the GaBP solution.  Let $R'_\gamma$ be defined as the edge-weight
matrix of a simple single-loop graph based on $\gamma$ with edge-weights defined
by
$r'_{\gamma_k,\gamma_{k+1}} = \frac{r_{\gamma_k,\gamma_{k+1}}}{1 -
\alpha_{\gamma_k \diff \gamma_{k+1}}}$. This construction is illustrated in
Figure 1(d).  Then,

\begin{lemma} $Z'_\gamma = \det (I-R'_\gamma)^{-1}
= \prod_{\ell \in \mathcal{L}_\gamma} Z_{\ell}$.
\end{lemma}
\emph{Proof.} Using Lemma 3, we see that the orbit-product $\prod_{\ell \in
\mathcal{L}_\gamma} Z_\ell$ is equal to the product over all non-trivial orbits
of the graph $G_\gamma$, that is, the product over all orbits in $G_\gamma$
which intersect the subgraph corresponding to $\gamma$. Using Lemma 1, this is
equivalent to computing the determinant of the corresponding submatrix of
$K_{G_\ell} = (I-R_{G_\ell})^{-1}$ where $R_{G_\ell}$ is the edge-weight matrix
of the computation graph.  This is equivalent to first eliminating each
computation tree (by Gaussian elimination/GaBP) attached to each node of
$\gamma$ to obtain a reduced graphical model $I-R_\gamma$ and then computing
$\det (I-R_\gamma)^{-1}$.  Using the GaBP solution, the effect of eliminating
each computation tree is to add a ``self-loop'' (diagonal element) to $R_\gamma$
with edge-weight $\alpha_{\gamma_k \diff \gamma_{k+1}} = \sum_{v \neq
\gamma_{k+1}} \alpha_{v,\gamma_k}$, obtained by summing the incoming messages to
node $\gamma_k$ from each of its neighbors in the subtree $T_{\gamma_k \diff
\gamma_{k+1}}$.  This elimination step is illustrated if Figure 1(b,c).  We may
use the orbit-product formula to compute the determinant. However, there are
infinitely many orbits in this graph due to the presence of a self-loop at each
of the remaining nodes.  At each node, an orbit may execute any number of steps
$m$ around this self-loop each with edge-weight $\alpha_{\gamma_k \diff
\gamma_{k+1}}$. Summing these, we obtain $\sum_{m=0}^\infty \alpha_{\gamma_k
\diff \gamma_{k+1}}^m = (1-\alpha_{\gamma_k \diff \gamma_{k+1}})^{-1}$.   
Hence, we can delete each self-loop and multiply the following edge's weight by
$(1-\alpha_{\gamma_k \diff \gamma_{k+1}})^{-1}$ and this preserves the value of
the determinant. This final reduction step is illustrated in Figure 1(c,d).
Then, the orbit-product $\prod_{\ell \in \mathcal{L}_\gamma} Z_\ell$ is equal to
$\det (I-R'_\gamma)^{-1}$ (e.g., based on the graph seen in Figure 1(d)). It is
straight-forward to compute the resulting determinant with respect to the single
(directed) cycle graph with edge-weights $R'_\gamma$.  There is only one orbit
in this graph and hence $\det (I-R'_\gamma)^{-1} = Z'_\gamma \triangleq 
(1-(R')^\gamma)^{-1}$ where $(R')^\gamma = \prod_{(ij) \in \gamma} 
(r'_{ij})^{n_{ij}(\gamma)}$ and $r'_{ij} = r_{ij} (1-\alpha_{i \diff j})^{-1}$.
\qed

\emph{Proof of Theorem 4.} Using these results, it is now simple
to show $\frac{Z}{Z^\bp}
= \prod_{\ell \not\in \mathcal{L}_\emptyset} Z_\ell
= \prod_{\gamma \neq \emptyset} \prod_{\ell \in \mathcal{L}_\gamma} Z_\ell
= \prod_{\gamma \neq \emptyset} \det(I-R'_\gamma)^{-1}
= \prod_{\gamma \neq \emptyset} (1-(R')^\gamma)^{-1}$. \qed

\section{Backtrackless Determinant Correction}

Next, we show that the correction factor \begin{displaymath} \frac{Z}{Z^\bp} =
\prod_{\gamma \neq \emptyset} Z'_\gamma = \prod_{\gamma \neq \emptyset}
\det(I-R'_\gamma)^{-1} \end{displaymath} may also be calculated as a
\emph{single determinant} based on the following \emph{backtrackless adjacency
matrix} of the graph.  We define $R' \in \mathbb{R}^{2|G| \times 2|G|}$ as
follows.  Let the rows and columns of $R'$ be indexed by directed edges $(ij)$
of the graph $G$.  Then, the elements of $R'$ are defined \begin{equation}
R'_{(ij),(kl)} = \left\{ \begin{array}{ll} r'_{kl}, & j=k \mbox{ and } i \neq
l\\ 0, & \mbox{otherwise}. \end{array} \right. \end{equation} This construction
is illustrated in Figure 2. Note that the walks generated by taking powers $R'$
correspond to backtrackless walks of the graph $G$.  The weight of an edge
$((ij)(jk))$ in $R'$ is defined as the (modified) edge-weight $r'_{jk}$ of the
endpoint $(jk)$.  The weight of an orbit in $R'$ may then be equivalently
defined as the product of \emph{node} weights $r'_{ij}$ taken over the orbit in
$R'$, which is equal to the weight of the corresponding backtrackless orbit of
$G$ (using the modified edge weights $r'_{ij}$).

\begin{figure}
(a) \hspace{.2cm} \input{grid.pstex_t} \hspace{.3cm} (b)
\hspace{-.2cm}\input{grid2.pstex_t}
\vspace{-.2cm}
\caption{(a) $3 \times 3$ grid $G$. (b) Graph $G'$ representing
the backtrackless adjacency matrix $R'$. Each node $ij$ represents
a directed edge of $G$, directed edges are drawn between nodes $ij$ and $jk$
which are
non-backtracking ($k \neq i$).}
\vspace{-.3cm}
\end{figure}

\begin{theorem} $Z = Z^\bp \times Z'$ where $Z' \triangleq \det(I-R')^{-1}$,
that is,
$\det(I-R)^{-1} = Z^\bp \times \det(I-R')^{-1}$.
\end{theorem}

Before providing the proof, we establish that walk-summability with respect to
$R$ implies walk-summability with respect to $R'$:

\begin{lemma} If $\rho(|R|)<1$ then $\rho(|R'|) \le \rho(|R|)$.  \end{lemma}

\emph{Proof.} Once the $\alpha$-parameters converge, the $\beta$-parameters
follow a linear system $\beta_{k+1} = R' \beta_k + b$ \cite{Moallemi09}. Hence,
the asymptotic convergence rate of GaBP is $\rho(R')$.  Compare this to the
Gauss-Jacobi (GJ) iteration $\mu_{k+1} = \mu_k + (h-J \mu_k) =
\sum_{t=0}^{k+1} R^t h$ ($\mu_0 = 0)$, which has convergence
rate $\rho(R)$.  It is clear that the GaBP iteration captures a superset of
those walks computed by GJ at each iteration (because the depth-$k$ computation
tree includes all $k$-length walks). Hence, for non-negative models ($R \ge 0$
and $h \ge 0)$ it must hold that the error in the GaBP estimate of $\mu$ is less
than or equal to the error of GJ (at every iteration). This implies $\rho(R')
\le \rho(R)$ if $R \ge 0$, from which we conclude $\rho(|R'|) \le \rho(|R|)$ in
walk-summable models.  \qed

\emph{Proof of Theorem 5.}  By construction, there is a one-to-one
weight-preserving correspondence between orbits of $G'$ and 
backtrackless orbits of $G$. The result then follows from the orbit-product 
representation of $Z'$ over $G'$  (Theorem 1, Lemma 5), which is 
equivalent to the backtrackless orbit-product of Theorem 4. \qed

One useful consequence of this result is that the error bound of Theorem 3 can
be improved to $\frac{1}{n} \left| \log \frac{Z^\bp}{Z} \right| = \frac{1}{n}
\left| \log Z' \right| \le \frac{\rho(|R'|)^g}{g(1-\rho(|R'|))}$.

It is impractical to compute the complete correction factor $Z' = \det
(I-R')^{-1}$, as this is not easier than directly computing $Z = \det
(I-R)^{-1}$.  However, because $R'$ is itself walk-summable, we can use this
representation as a starting point for constructing approximate corrections
such as the one considered in the next section.

\section{Block-Resummation Method}

Next, we consider an efficient method to approximate $Z(A)=\det(I-A)^{-1}$ for
walk-summable models $\rho(|A|)<1$.  This method can be used to either directly
approximate $Z(R)$ ($A=R$) or to approximate the GaBP-correction $Z'$ ($A=R')$.

Given a graph $G$ based on vertices $V$, we specify a set of blocks $\mathcal{B}
= (B_k \subset V, k=1,\dots,|\mathcal{B}|)$ chosen such that: (1) Every short
orbit $|\ell| < L$ is covered by some block $B \in \mathcal{B}$, and (2) If
$B,B' \in \mathcal{B}$ then $B \cap B' \in \mathcal{B}$. We also define block
weights $w_B$ as follows:  $w_B=1$ for maximal blocks (not contained by another
block) and $w_B = 1 - \sum_{B' \supsetneq B} w_{B'}$ for
non-maximal blocks (these weights may be negative). This insures that
$\sum_{B' \supseteq B} w_{B'} = 1$ for each $B \in \mathcal{B}$. Then, we define
our estimate
\begin{equation}
Z_{\mathcal{B}} \triangleq \prod_B Z_B^{w_B} \triangleq \prod_B
(\det(I-A_B)^{-1})^{w_B}
\end{equation}
where $A_B$ denotes the $|B| \times |B|$ principle submatrix of $A$
corresponding to $B$. 

This approximation method is similar in spirit to approximations used
elsewhere (e.g., Kikuchi approximations to free-energy \cite{Yedidia05}). 
However, the new insights offered by the orbit-product view allows us to
give our estimate a precise interpretation in walk-summable Gaussian 
models:

\begin{theorem} $Z_{\mathcal{B}} = \prod_{\ell \in \mathcal{L}_{\mathcal{B}}}
Z_\ell$
where $\mathcal{L}_{\mathcal{B}} \triangleq \cup_{B \in \mathcal{B}}
\mathcal{L}_B$ and
$\mathcal{L}_B$ is the set of all orbits covered by $B$.
\end{theorem}
\emph{Proof.} $Z_{\mathcal{B}} = \prod_B \prod_{\ell \in B} Z_\ell^{w_B} =
\prod_{\ell \in \mathcal{L}_{\mathcal{B}}} Z_\ell^{\sum_{B \supset \ell} w_B} =
\prod_{\ell \in \mathcal{L}_{\mathcal{B}}} Z_\ell$ where $\sum_{B \supset \ell}
w_B = 1$ follows from the definition of the block weights. \ignore{Let $B'$
denote the smallest block that covers $\ell$.  Then, $\sum_{B \supset \ell} w_B
= \sum_{B \subseteq B'} w_B = 1$.} $\diamond$

Moreover, we can then bound the error of the estimate.  Noting
that $\mathcal{L}_\mathcal{B}$ includes all short orbits of the
graph, we can derive the following result by a similar proof as 
for Theorem 3:

\begin{corollary} $\tfrac{1}{n} \left| \log \frac{Z_{\mathcal{B}}}{Z} \right|
\le
  \frac{\rho(|A|)^L}{L(1-\rho(|A|))}$.
\end{corollary}

Thus, for the class of models with $\rho < 1$, we obtain an approximation scheme
which converges to the correct determinant as the parameter $L$ is made large
with error decaying exponentially in $L$. The estimate $Z_{\mathcal{B}}(R)$
includes all orbits that are covered by some block.  The improved GaBP-based
estimate $Z^\bp Z_{\mathcal{B}}(R')$ includes all orbits $\ell$ such that
$\gamma = \Gamma(\ell)$ is covered by some block.  Thus, the GaBP-based
correction includes many more orbits. We also note that the error-bound using
the GaBP-based estimate is typically smaller as we have shown that $\rho(|R'|)
\le \rho(|R|)$ (if $\rho(|R|)<1$).

\begin{figure}
\centering
\hspace{.16cm}\epsfig{file=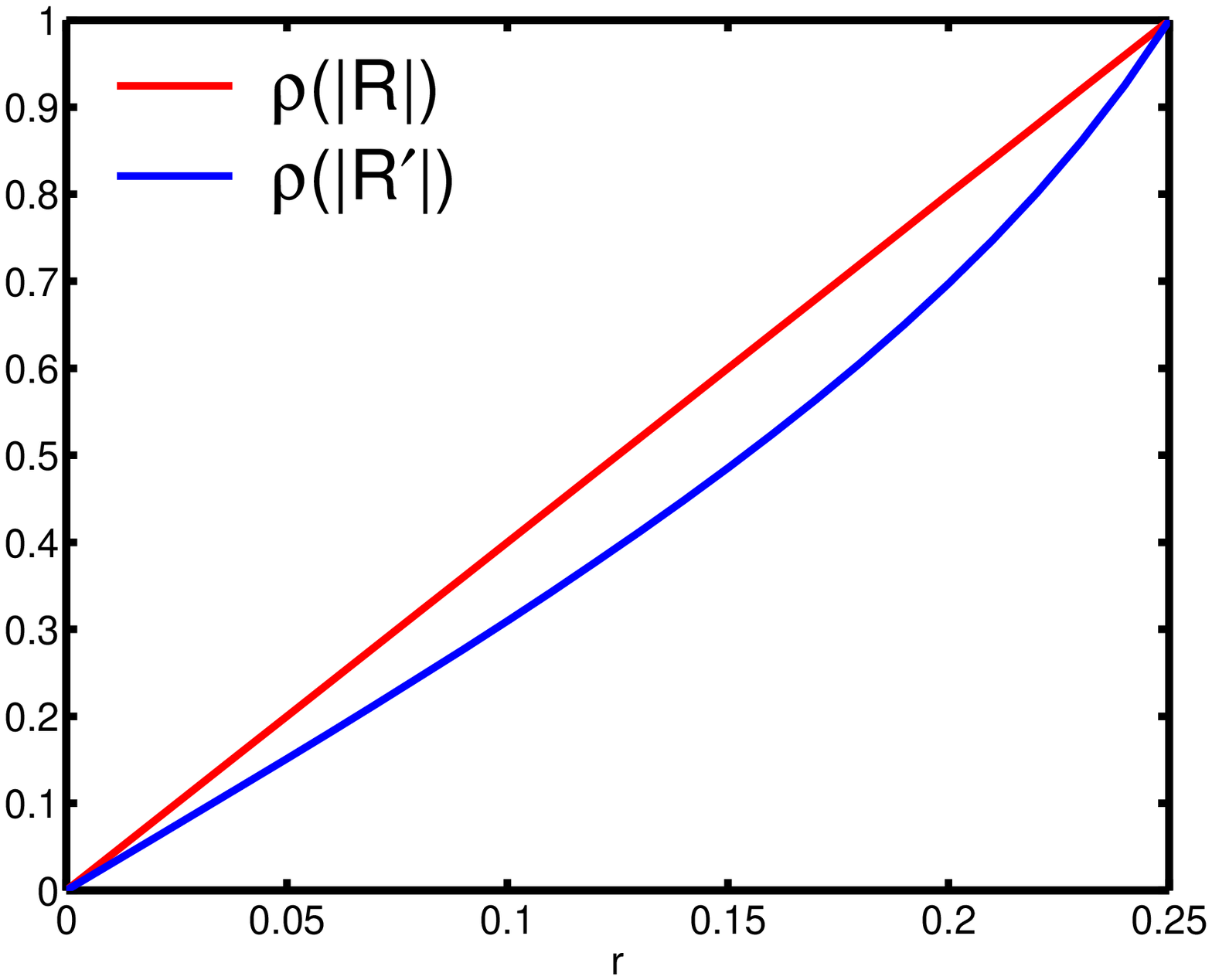,scale=0.2}\hspace{.2cm}
\epsfig{file=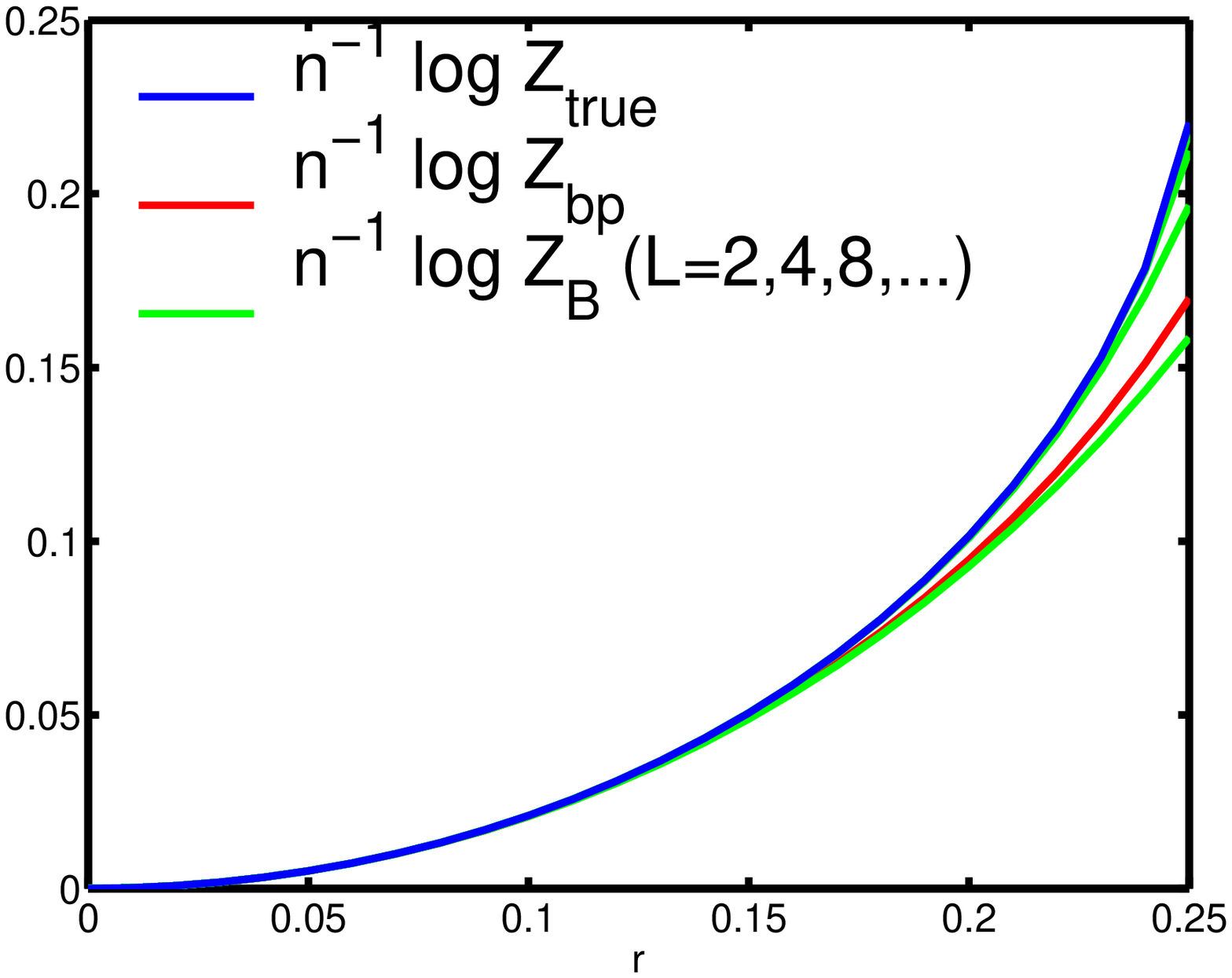,scale=0.2}\\
\vspace{-.1cm}
\hspace{.2cm}(a)\hspace{3.6cm}(b)\\
\epsfig{file=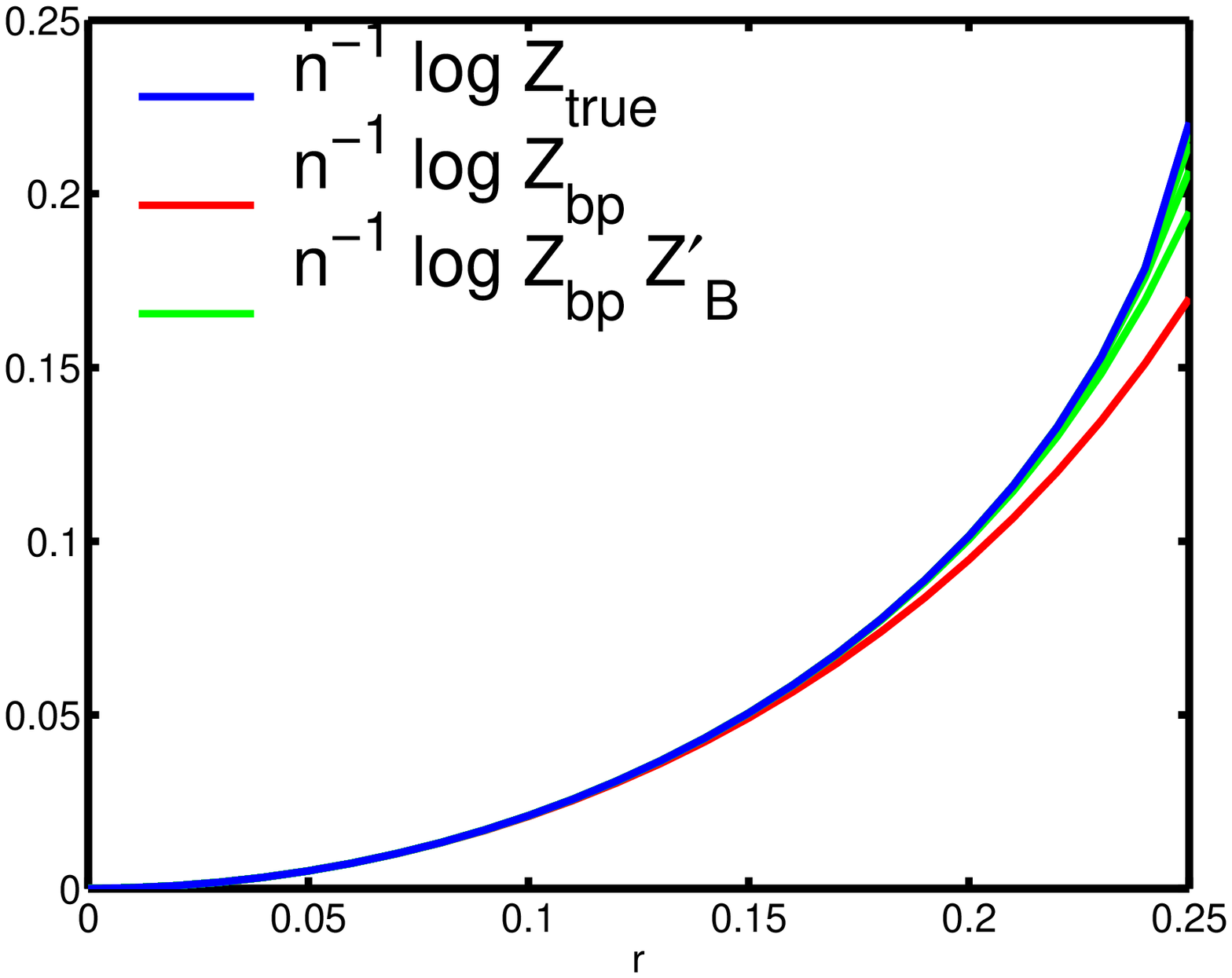,scale=0.2}\hspace{.2cm}
\epsfig{file=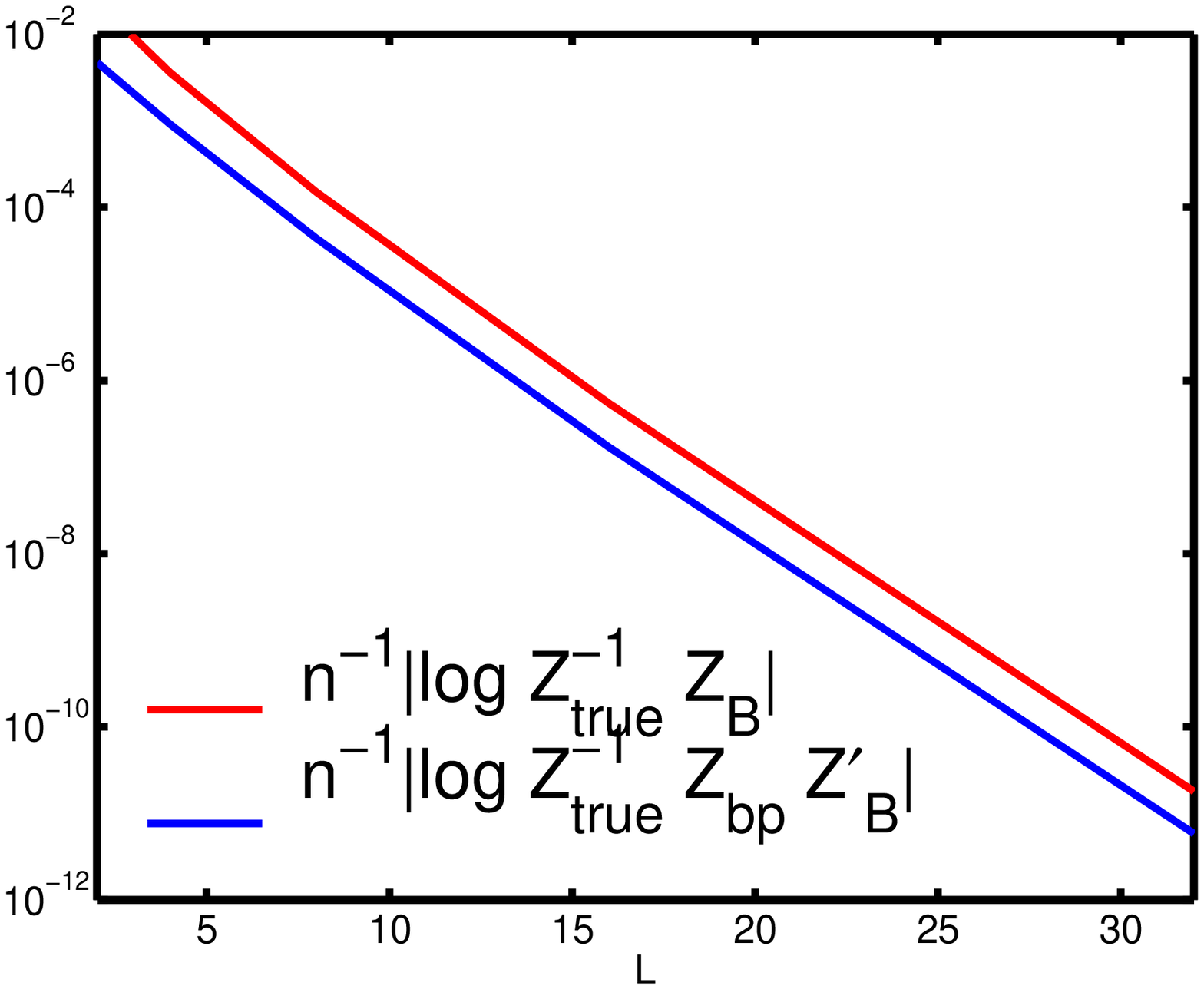,scale=0.2}\\
\vspace{-.1cm}
\hspace{.2cm}(c)\hspace{3.6cm}(d)\\
\vspace{-.2cm}
\caption{Demonstration of determinant approximation method for $256\times256$
periodic grid with uniform edge weights $r \in (0,.25)$. Plots of (a)
$\rho(|R|)$ and $\rho(|R'|)$ vs $r$; (b) ($\tfrac{1}{n} \log$ of) $Z$, $Z^\bp$
and $Z_\mathcal{B}$ (with $L = 2,4,8,16,32$) vs $r$; (c) $\tfrac{1}{n} \log
(Z^\bp Z'_\mathcal{B})$ vs $r$; and (d) $\tfrac{1}{n} |\log
(Z^{-1}Z_\mathcal{B})|$ and $\tfrac{1}{n}|\log (Z^{-1}Z^\bp Z'_\mathcal{B})|$ vs
$L=2,4,8,16,32$ for $r=.23$. In (b) and (c) the estimates for $L=8,16,32$ are
all nearly exact (and therefore hard to distinguish in the plot) and the errors
are largest near the walk-summable threshold $r=.25$. Estimates (b) are not
based on GaBP and are actually worst than $Z^\bp$ for $L=2$. However, the
GaBP-corrections (c) are strictly better than $Z^\bp$.} \vspace{-.3cm}
\end{figure}

\paragraph{Construction of $\mathcal{B}$ for Grids}

To achieve an error bound $\tfrac{1}{n} |\log\tfrac{Z_{\mathcal{B}}}{Z}| \le
\varepsilon$ we must choose $L \sim \log \varepsilon^{-1}$.  Then, the
computation needed to achieve this precision will depend on both the number of
blocks and the block size needed to cover all orbits up to this length.  In
certain classes of sparse graphs, it should be possible to control the
complexity of the method.  As an example, we demonstrate how to choose blocks
for $2D$ grids. Consider the $\sqrt{n} \times \sqrt{n}$ square grid in which
each vertex is connected to its four nearest neighbors.  We may cover this graph
by $L \times L$ blocks shifted (both vertically and horizontally) in increments
of $\tfrac{L}{2}$ (let $L$ be even).  It can be seen that this set of blocks
covers all loops shorter than $L$. To include all intersections of blocks, we
add $L \times \tfrac{L}{2}$, $\tfrac{L}{2} \times L$ and $\tfrac{L}{2} \times
\tfrac{L}{2}$ blocks.  The block weights are $w_{L \times L} = 1$, $w_{L \times
L/2} = w_{L/2 \times L} = -1$ and $w_{L/2 \times L/2} = 1$.  The complexity of
computing the determinant of an $L \times L$ block is $\mathcal{O}(L^3)$ and the
total number of blocks is $\mathcal{O}(n/L^2)$.  Hence, the total complexity is
$\mathcal{O}(nL) = \mathcal{O}(n \log \varepsilon^{-1})$.  \ignore{Hence, in
grids we achieve both linear complexity in the number of variables and linear
complexity in the precision $\log \varepsilon^{-1}$ (which roughly corresponds
to the number of significant digits of our estimate $\tfrac{1}{n} \log Z$).}

We test our approach numerically on a $256\times256$ square grid (with periodic
boundary conditions).  We set all edge weights to $r$ and test the quality of
approximation using both estimates $Z_\mathcal{B}(R)$ and $Z^\bp
Z_\mathcal{B}(R')$ for $r \in (0,.25)$ ($J=I-R$ becomes indefinite for
larger values of $r$) and block sizes $L = 2,4,8,16,32$.  The results are shown
in Figure 3.  As expected, accuracy rapidly improves with increasing $L$ in both
methods and the GaBP-correction approach is more accurate.

\section{Conclusion and Future Work}

We have demonstrated an orbit-product representation of the determinant (the
partition function of the Gaussian model) and interpreted the estimate obtained
by GaBP as corresponding to totally backtracking orbits.  Furthermore, we have
shown how to correct the GaBP estimate in various ways which involve
incorporating backtrackless orbits (e.g. cycles) of the graph. In particular, we
demonstrated an efficient block-resummation method to compute truncated
orbit-products in sparse graphs (demonstrated for grids). These methods also
extend to address estimation of the matrix inverse (the covariance matrix of the
Gaussian model), which may in turn be used as an efficient preconditioner for
iterative solution of linear systems.  We leave these extensions for a longer
report.

In future work, we plan to extend the method of constructing an efficient set of
blocks to other classes of sparse graphs.   It may also be fruitful to extend
our analysis to generalized belief propagation \cite{Yedidia05} in Gaussian
models. In a related direction, we
intend to explore methods to ``bootstrap'' GaBP using the factorization $Z(R) =
\left(\prod_{k=0}^\infty Z(-R^{2^k}))\right)^{-1}$, which follows from the
formula $(I-R)^{-1} = \prod_k (I+R^{2^k})$. By computing $Z^\bp(-R^{2^k})$ for
small values of $k$ we may capture short backtrackless orbits of the graph.
Another direction is to investigate generalization of the formula $Z = Z^\bp Z'$
to non-walksummable models, perhaps using methods of \cite{Chernyak08}.  A
related idea is to approximate a non-walksummable model by a walk-summable one
and then correct estimates obtained from the walk-summable model to better
approximate the non-walksummable model.

\bibliography{jcc-icml09}

\begin{thebibliography}{10}

\bibitem{Bickson08b}
D.~Bickson, D.~Dolev, and E.~Yom-Tov.
\newblock A {G}aussian belief propagation solver for large scale {SVM}s.
\newblock In {\em 5th Europ. Conf. Complex Systems}, 2008.

\bibitem{Bickson08a}
D.~Bickson, O.~Shental, P.~Siegel, J.~Wolf, and D.~Dolev.
\newblock {G}aussian belief propagation based multiuser detection.
\newblock In {\em IEEE Int. Symp. on Inform. Th.}, pages 1878--1882, 2008.

\bibitem{Chernyak08}
V.~Chernyak and M.~Chertkov.
\newblock Fermions and loops on graphs {I}: Loop calculus for determinant.
\newblock {\em J. Statistical Mechanics: Theory and Experiments}, 2008.

\bibitem{Chertkov06}
M.~Chertkov and V.~Chernyak.
\newblock Loop series for discrete statistical models on graphs.
\newblock {\em J. Statistical Mechanics: Theory and Experiments}, 2006.

\bibitem{Cseke08}
B.~Cseke and T.~Heskes.
\newblock Bounds on the {B}ethe free energy for {G}aussian networks.
\newblock In {\em Uncertainty in Artificial Intelligence}, pages 97--104, 2008.

\bibitem{Johnson06}
J.~Johnson, D.~Malioutov, and A.~Willsky.
\newblock Walk-sum interpretation and analysis of {G}aussian belief
  propagation.
\newblock In {\em Adv. in Neural Inform. Processing Systems 18}, pages
  579--586, 2006.

\bibitem{Malioutov06}
D.~Malioutov, J.~Johnson, and A.~Willsky.
\newblock Walk-sums and belief propagation in {G}aussian graphical models.
\newblock {\em J. of Machine Learning Research}, 7:2031--2064, 2006.

\bibitem{Moallemi06}
C.~Moallemi and B.~{Van Roy}.
\newblock Consensus propagation.
\newblock {\em IEEE Trans. on Inform. Th.}, 52:4753--4766, 2006.

\bibitem{Moallemi09}
C.~Moallemi and B.~{Van Roy}.
\newblock Convergence of min-sum message passing for quadratic optimization.
\newblock {\em IEEE Trans. on Inform. Th.}, 55:2413--2423, May 2009.

\bibitem{Plarre04}
K.~Plarre and P.~Kumar.
\newblock Extended message passing algorithm for inference in loopy {G}aussian
  graphical models.
\newblock {\em Ad Hoc Networks}, 2:153--169, 2004.

\bibitem{Rusmevichientong01}
P.~Rusmevichientong and B.~{Van Roy}.
\newblock An analysis of belief propagation on the turbo decoding graph with
  {G}aussian densities.
\newblock {\em IEEE Trans. on Inform. Th.}, 47:745--765, 2001.

\bibitem{Stark96}
H.~Stark and A.~Terras.
\newblock Zeta functions of finite graphs and coverings.
\newblock {\em Adv. in Math.}, 121:124--165, 1996.

\bibitem{Weiss01}
Y.~Weiss and W.~Freeman.
\newblock Correctness of belief propagation in {G}aussian graphical models of
  arbitrary topology.
\newblock {\em Neural Computation}, 13:2173--2200, 2001.

\bibitem{Yedidia05}
J.~Yedidia, W.~Freeman, and Y.~Weiss.
\newblock Constructing free-energy approximations and generalized belief
  propagation algorithms.
\newblock {\em IEEE Trans. Inform. Th.}, 51:2282--2312, 2005.

\end{thebibliography}
\bibliographystyle{plain}

\end{document}